
\magnification1200
\input amstex.tex
\documentstyle{amsppt}

\vsize=18cm

\footline={\hss{\vbox to 2cm{\vfil\hbox{\rm\folio}}}\hss}
\nopagenumbers
\def\DJ{\leavevmode\setbox0=\hbox{D}\kern0pt\rlap
{\kern.04em\raise.188\ht0\hbox{-}}D}

\def\txt#1{{\textstyle{#1}}}
\baselineskip=13pt
\def\hf{{\textstyle{1\over2}}}
\def\a{\alpha}\def\b{\beta}
\def\d{{\,\roman d}}
\def\e{\varepsilon}
\def\f{\varphi}
\def\G{\Gamma}
\def\k{\kappa}
\def\s{\sigma}
\def\t{\theta}
\def\={\;=\;}

\def\zt{\zeta(\hf+it)}

\def\D{\Delta}

 \def\t{\theta}
\def\hf{{\textstyle{1\over2}}}
\def\txt#1{{\textstyle{#1}}}
\def\f{\varphi}

\font\tenmsb=msbm10
\font\sevenmsb=msbm7
\font\fivemsb=msbm5
\newfam\msbfam
\textfont\msbfam=\tenmsb
\scriptfont\msbfam=\sevenmsb
\scriptscriptfont\msbfam=\fivemsb
\def\Bbb#1{{\fam\msbfam #1}}

\def \NN {\Bbb N}

\def \ZZ {\Bbb Z}

\font\ff=cmr8
\def\txt#1{{\textstyle{#1}}}
\baselineskip=13pt

\font\teneufm=eufm10
\font\seveneufm=eufm7
\font\fiveeufm=eufm5
\newfam\eufmfam
\textfont\eufmfam=\teneufm
\scriptfont\eufmfam=\seveneufm
\scriptscriptfont\eufmfam=\fiveeufm
\def\mathfrak#1{{\fam\eufmfam\relax#1}}

\font\tenmsb=msbm10
\font\sevenmsb=msbm7
\font\fivemsb=msbm5
\newfam\msbfam
     \textfont\msbfam=\tenmsb
      \scriptfont\msbfam=\sevenmsb
      \scriptscriptfont\msbfam=\fivemsb
\def\Bbb#1{{\fam\msbfam #1}}

\def \NN {\Bbb N}

\def \ZZ {\Bbb Z}

  \def\rightheadline{{\hfil{\ff
  On the higher moments of the error term in the divisor problem}
\hfil\tenrm\folio}}

  \def\leftheadline{{\tenrm\folio\hfil{\ff
   A. Ivi\'c and P. Sargos}\hfil}}
  \def\emptyheadline{\hfil}
  \headline{\ifnum\pageno=1 \emptyheadline\else
  \ifodd\pageno \rightheadline \else \leftheadline\fi\fi}

\topmatter
\title
ON THE HIGHER MOMENTS OF THE ERROR TERM IN THE DIVISOR PROBLEM
\endtitle
\author   Aleksandar Ivi\'c and Patrick Sargos \endauthor
\address
Aleksandar Ivi\'c, Katedra Matematike RGF-a
Universiteta u Beogradu, \DJ u\v sina 7, 11000 Beograd,
Serbia (Yugoslavia).
\medskip
Patrick Sargos, Institut Elie Cartan,
Universit\'e Henri Poincar\'e, BP. 239,
54506 Nancy, France
\bigskip
\endaddress
\keywords
Dirichlet divisor problem, mean fourth power, asymptotic formula
\endkeywords
\subjclass
11N37, 11M06 \endsubjclass
\email {\tt
aivic\@rgf.bg.ac.yu,  sargos\@iecn.u-nancy.fr} \endemail
\dedicatory
\enddedicatory
\abstract
{Let $\D(x)$ denote the error term in the Dirichlet
divisor problem. Our main results are the asymptotic formulas
$$
\int_1^X\D^3(x)\d x = BX^{7/4} + O_\e(X^{\b+\e})  \qquad(B > 0)
$$
and
$$
\int_1^X\D^4(x)\d x = CX^2 + O_\e(X^{\gamma+\e})  \qquad(C > 0)
$$
with $\b = {7\over5}, \gamma = {23\over12}$. This improves on the values
$\b = {47\over28}, \gamma  = {45\over23}$, due to  K.-M. Tsang.
A result on integrals of $\D^3(x)$ and $\D^4(x)$ in short
intervals is also proved.}
\endabstract
\endtopmatter

\head
1. Introduction and statement of results
\endhead
For a fixed $k\in\NN$, let
$$
\D_k(x) = \sum_{n\le x}d_k(n) - xP_{k-1}(\log x)\leqno(1.1)
$$
denote the error term in the (general) Dirichlet divisor
problem. Here $d_k(n)$ denotes the number of ways
$n$ may be written as a product of $k$ factors
(so that $d_2(n) = d(n)$ is the number of divisors
of $n$), and  $P_{k-1}(z)$ is a suitable polynomial of
degree $k-1$ in $z$ (see e.g., [6, Chapter 13] for more
details). In particular,
$$
\D_2(x) \equiv \D(x) =   \sum_{n\le x}d(n) - x(\log x + 2\gamma - 1)
\leqno(1.2)
$$
is the error term in the classical Dirichlet divisor problem
($\gamma = -\G'(1) = 0.5772\ldots\,$ is Euler's constant).
A vast literature exists on the estimation of $\D_k(x)$ and
especially on $\D(x)$ (see e.g., [6] and [22]),
both pointwise and in various means.
Here we shall be concerned with the third and fourth moment of
$\D(x)$.
We note that the first author in [5] proved
a large values estimate for $\D(x)$, which yielded the bound
$$
\int_1^X\D^4(x)\d x \ll_\e X^{2+\e},\leqno(1.3)
$$
where here and later $\e$ denotes arbitrarily small, positive
constants, which are not necessarily the same ones at each
occurrence. The asymptotic formula for the fourth moment
with an error term was obtained by
K.-M. Tsang [23]. He proved that
$$
\int_1^X\D^4(x)\d x = CX^2 + O_\e(X^{\gamma+\e})  \leqno(1.4)
$$
holds with explicitly given $C\,(>0)$ and $\gamma = 45/23 = 1.956\ldots\;$.
Tsang also
proved an asymptotic formula for the integral of the cube of $\D(x)$,
namely
$$
\int_1^X\D^3(x)\d x = BX^{7/4} + O_\e(X^{\b+\e})  \leqno(1.5)
$$
with explicit $B>0$ and $\b = 47/28 = 1.678\ldots\;$.
Later D.R. Heath-Brown [3]  established asymptotic formulas
for higher moments of $\D(x)$, but his method does not produce
error terms.
\smallskip
The main  aim of this paper is to improve on the values
of Tsang's exponents $\b$ and $\gamma$ in (1.5) and (1.4), respectively.
The results are

\bigskip
THEOREM 1. {\it We have}
$$
\int_1^X\D^3(x)\d x = BX^{7/4} + O_\e(X^{\b+\e})
 \qquad\left(\b = {7\over5} = 1.4\,,
\;B > 0\right).\leqno(1.6)
$$

\bigskip
THEOREM 2. {\it We have}
$$
\int_1^X\D^4(x)\d x = CX^2 + O_\e(X^{\gamma+\e})
\qquad\left(\gamma = {23\over12} = 1.91666\ldots\,,
\;C > 0\right).\leqno(1.7)
$$

\medskip
The true values of $\b$ and $\gamma$ for which (1.6) and (1.7) hold
are hard to determine. However, it is not difficult to show that
(1.7) cannot hold with $\gamma < 5/4$. To see this, note first
that, for $1 \ll H \ll X$,
$$\eqalign{\cr&
\D(X) - {1\over 2H}\int_{X-H}^{X+H}\D(x)\d x\cr&
= {1\over 2H}\int_{X-H}^{X+H}(\D(X) - \D(x))\d x\cr&
\ll {1\over H}\int_{X-H}^{X+H}\Bigl(\bigl|\sum_{X\le n\le x}d(n)\bigr|
+ O(H\log X)\Bigr)\d x
\cr& \ll_\e HX^\e,\cr}
$$
which by H\"older's inequality and (1.7) gives
$$\eqalign{\cr&
\D^4(X) \ll_\e {1\over H}\int_{X-H}^{X+H}\D^4(x)\d x + H^4X^\e \cr&
\ll_\e X + H^{-1}X^{\gamma+\e} + H^4X^\e.\cr}
$$
But taking $H = X^{\gamma/5}$ we get
$$
\D(X) \ll X^{1/4} + X^{\gamma/5+\e},
$$
which is a contradiction if $\gamma < 5/4$, since $\limsup_{X\to\infty}
|\D(X)|X^{-1/4}=\infty$ is a classical result of G.H. Hardy
[2]. For the cube
this procedure does not work directly, since $\D^3(x)$ may be negative.

\medskip
We note that asymptotic formulas for moments of $\D(x)$ in short
intervals were recently investigated by W.G. Nowak
[20], but his results do not imply the asymptotic formulas (1.4) and (1.5).
Nowak actually works, for technical reasons, with $\D(x^2)$ instead of
$\D(x)$.  With a slight change of  notation, his formulas
for the cube and the fourth moment may be written as
$$
\int_{X-H}^{X+H}\D^3(x)\d x = (D_1 + o(1))HX^{3/4}
\qquad(X^{3/4+\delta} \le H \le \lambda X,\, 0 < \lambda \le 1)
\leqno(1.8)
$$
for any fixed $0 < \delta < {1\over4}$ and $X\to \infty$, and similarly
$$
\int_{X-H}^{X+H}\D^4(x)\d x = (D_2 + o(1))HX
\qquad(X^{3/4+\delta} \le H \le \lambda X,\, 0 < \lambda \le 1),
\leqno(1.9)
$$
where $D_j = D_j(\lambda)\;(>0)$ is explicitly given by Nowak
for $j = 1,2$. For example, $D_1 = 7B/2$ if $H = o(X)$ and
$D_1 = B\lambda^{-1}((1+\lambda)^{7/4} - (1-\lambda)^{7/4})$
if $H = \lambda X$, where $B$ is the constant appearing in (1.6).
We shall improve the range for which (1.8) and
(1.9) hold. We avoid $``o(1)"$ in (1.9) and formulate our results
as

\bigskip
THEOREM 3. {\it For any fixed $0 < \delta < {1\over3}$ there exists
$\kappa > 0$ such that uniformly}
$$
\int_{X}^{X+H}\D^3(x)\d x = B\left((X+H)^{7/4} - X^{7/4}\right)
(1 + O(X^{-\kappa}))
\quad(X^{7/12+\delta} \le H \le  X),
\leqno(1.10)
$$
{\it and also}
$$
\int_{X}^{X+H}\D^4(x)\d x = C\left((X+H)^2-X^2\right)(1 + O(X^{-\kappa}))
\quad(X^{2/3+\delta} \le H \le  X),
\leqno(1.11)
$$
{\it where B and C are the constants appearing in Theorem} 1 {\
\it and Theorem} 2.

\medskip
From the proof of Theorem 3 it will be clear that, using the sharpest
bound for $\D(x)$ and estimates for certain exponentilal sums,
we can improve on the exponents 7/12 and 2/3 which
appear in (1.10) and (1.11), respectively.
Following the method of proof of Theorem 1 and Theorem 2, one can obtain
analogous results for the cube and fourth power of two well-known
number theoretic error terms. This is given by
\bigskip
{\bf Corollary 1}. {\it We have}
$$\eqalign{\cr
\int_1^TE^3(t)\d t &= B_1T^{7/4} + O_\e(T^{\b_1+\e}),\cr
\int_1^XP^3(x)\d x &= B_2X^{7/4} + O_\e(X^{\b_2 +\e}),\cr}
\leqno(1.12)
$$
{\it with       }$\b_1 = 5/3,  \b_2 = 7/5$.
\bigskip
{\bf Corollary 2}. {\it We have}
$$\eqalign{\cr
\int_1^TE^4(t)\d t &= C_1T^2 + O_\e(T^{\gamma_1+\e}),\cr
\int_1^XP^4(x)\d x &= C_2X^2 + O_\e(X^{\gamma_2 +\e}),\cr}
\leqno(1.13)
$$
{\it with }$ \gamma_1 = \gamma_2 = 23/12$.

\bigskip
Here $B_1,B_2,C_1,C_2$ are explicit, positive constants,
$$
E(T) \;=\;\int_0^T|\zt|^2\d t - T\left(\log\left({T\over2\pi}\right)
 + 2\gamma - 1\right)
$$
denotes the error term in the mean square formula for $|\zt|$, while
$$
P(x) = \sum_{n\le x}r(n) - \pi x\qquad(r(n) = \sum_{n=a^2+b^2;a,b\in\ZZ}1)
$$
denotes the error term in the circle problem.

Namely we have the explicit, truncated formula (see e.g., [6] or [22])
$$
\D(x) = {1\over\pi\sqrt{2}}
x^{1\over4}\sum_{n\le N}d(n)n^{-{3\over4}}\cos(4\pi\sqrt{nx}
- {\txt{1\over4}}\pi) +
O_\e(x^{{1\over2}+\e}N^{-{1\over2}})\quad(2 \le N \ll x).
\leqno(1.14)
$$
All our results on $\D(x)$ will be proved by using solely (1.14), without
recourse to the arithmetic structure of $d(n)$, except the trivial
bound $d(n) \ll_\e n^\e$. For $P(x)$ there also exists an explicit
formula, namely
$$
P(x) = -{1\over\pi}x^{1\over4}\sum_{n\le N}r(n)n^{-{3\over4}}
\cos(2\pi\sqrt{nx} +
{\txt{1\over4}}\pi) + O_\e(x^{{1\over2}+\e}N^{-{1\over2}})\;(2 \le N \ll x).
\leqno (1.15)
$$
This is analogous to (1.14), and the proof
follows on the same lines, e.g., by the method given in Titchmarsh [22].
Hence all the results on $\D(x)$ that depend only on (1.14) and
$d(n) \ll_\e n^\e$ have their corresponding analogues for $P(x)$, so that
 by this principle the second formulas in (1.12) and (1.13) follow.

\smallskip
In what concerns the formulas involving $E(t)$, we note that
$$\eqalign{\cr
\int_0^TE^3(t)\d t &\;=\; 16\pi^4\int_0^{T\over2\pi}(\D^*(t))^3\d t
+ O(T^{5/3}\log^{3/2}T),\cr
\int_0^TE^4(t)\d t &\;=\; 32\pi^5\int_0^{T\over2\pi}(\D^*(t))^4\d t
+ O(T^{23/12}\log^{3/2}T)\cr}\leqno(1.16)
$$
holds, as proved by the first author [7, Theorem 2]. In (1.16) we set
$$
\D^*(x) := - \D(x) + 2\D(2x) - \hf\D(4x).
$$
Then  the arithmetic interpretation of $\D^*(x)$ is
$$
\hf\sum_{n\le4x}(-1)^nd(n) \;=\; x(\log x + 2\gamma - 1) + \D^*(x).
$$
One also has (see  [6, eq. (15.68)]), for $1 \ll N \ll x$,
$$
\D^*(x) = {1\over\pi\sqrt{2}}x^{1/4}\sum_{n\le N}(-1)^nd(n)n^{-3/4}
\cos(4\pi\sqrt{nx} - {\txt{1\over4}}\pi) + O_\e(x^{1/2+\e}N^{-1/2}),
$$
which is completely analogous to (1.14). Hence the analogues of our
formulas (1.4) and (1.5) hold for $\D^*(x)$, and therefore by
(1.14) we easily obtain then the first formulas in (1.10) and (1.11).
Note that the exponents $\b_1$ and $\b_2$ in (1.12) are not equal;
this reflects the state of art that in (1.16) the first error term
has the exponent of $T$ equal to 5/3, and $7/5 < 5/3$.

\smallskip
It should be remarked that moments of other number-theoretic
error terms can be dealt with by our methods. This involves,
for example, $A(x) = \sum_{n\le  x}a(n)$, where $a(n)$ is the
$n$-th Fourier coefficient of a holomorphic cusp form of
weight $\kappa = 2n \,(\ge 12)$, which possesses an explicit formula
analogous to (1.14). In general, error terms connected with
the coefficients of Dirichlet series belonging to the well-known
Selberg class of degree two may be considered.

\medskip
The plan of the paper is as follows.  In Section 2 we shall prove
lemmas on the spacing of square roots, which are needed for the proof
of our results.  In Section 3 we shall prove
Theorem 1. Section 4 contains the proof of Theorem 2,
while Theorem 3 will be proved in Section 5.

\medskip
{\bf Remarks and acknowledgements}. After the first version of this
paper was written, with the exponent 7/5 in (1.6), Kai-Man Tsang
kindly informed us of the doctoral thesis of his student Yuk-Kam Lau [8].
Lau investigates a slightly more general function than $\D(x)$, namely
$\D_a(x)$, the error term in the asymptotic formula for the summatory
function of $\s_a(n) = \sum_{d|n}d^a$ for $a$ in a certain range.
Lau obtains, for the integral of $\D_a^3(x)$ in the case $a = 0$, the
same exponent $7/5+\e$ as we do, ``following the method of Tsang and some
refinements suggested by him" ([8, p. 98]).
Since Lau's result has not been published in any periodical,
and we obtained (1.6) independently of him, we thought appropriate to
retain our proof of (1.6) in the final version of the paper.

We are also grateful to T. Meurman who  informed us of
the papers of W. Zhai [24], who kindly sent us his works and
made valuable remarks.
Zhai establishes asymptotic formulas for the integrals of $\D^j(x)$,
when $3 \le j \le 9$. For $j = 3,4$ in the notation of (1.6) and (1.7)
he had $\b = 3/2$ and $\gamma = 80/41$, which is poorer than what
we obtained, although in correspondence Zhai indicated that his methods
also can yield $\b = 7/5$. Finally we thank W.G. Nowak for sending us
a preprint of [20].

\head
2. Lemmas on the spacing of the square roots
\endhead

Both in the proof of Theorem 1 and Theorem 2 the basic approach is
obvious: the sum approximating $\D(x)$ (cf. (1.12)) is raised to the
 third, respectively fourth power, and the resulting expressions are
integrated. In this process sums and differences of square roots
will appear in the exponentials. Thus several lemmas on the
spacing of the square roots will be needed. We begin with

\medskip
LEMMA 1 (O. Robert--P. Sargos [21]). {\it Let $k\ge 2$ be a fixed
integer and $\delta > 0$ be given.
Then the number of integers $n_1,n_2,n_3,n_4$ such that
$N < n_1,n_2,n_3,n_4 \le 2N$ and}
$$
|n_1^{1/k} + n_2^{1/k} - n_3^{1/k} - n_4^{1/k}| < \delta N^{1/k}
$$
{\it is, for any given $\e>0$,}
$$
\ll_\e N^\e(N^4\delta + N^2).\leqno(2.1)
$$

\medskip
LEMMA 2. {\it If m,n,k are natural numbers such that $\sqrt{m}+\sqrt{n}
 \not = \sqrt{k}$, then}
$$
|\sqrt{m}+\sqrt{n} - \sqrt{k}| \gg (mnk)^{-1/2}.\leqno(2.2)
$$
{\it If m,n,l,k are natural numbers such that $m\le k,\; n\le k$ and
 $\sqrt{m}+\sqrt{n} \pm \sqrt{k} \not= \sqrt{l}$, then}
$$
|\sqrt{m}+\sqrt{n} \pm \sqrt{k} - \sqrt{l}| \gg  k^{-2}(mnl)^{-1/2}.
\leqno(2.3)
$$

\medskip {\bf Proof.} These  results
should be compared with Tsang [23, Lemma 2] and [23, Lemma 3], who had
(2.2) and (2.3) with the right-hand sides replaced by $\max(m,n,k)^{-1/2}$
and $\max(m,n,k,l)^{-7/2}$, respectively.
Thus our bounds are better when at least one
of the integers in question is smaller than the other ones.

To prove (2.2), we note first that if $4x\in \NN$ is not a square, then
$$
||2\sqrt{x}|| \;\gg\;{1\over{\sqrt{x}}},\leqno(2.4)
$$
where as usual $||y||$ is the distance of $y$ to the nearest integer.
Namely if $\eta = ||2\sqrt{x}|| $, then $2\sqrt{x} = n \pm \eta,\,n\in\NN$
and $0 < \eta \le \hf$. Then $4x = n^2 \pm 2\eta n + \eta^2$.
If $4x$ is not a square,  we have
$$
1 \le |4x - n^2| \le 2\eta n + \eta^2 \ll \eta\sqrt{x},
$$
which implies (2.4).

Now let $\t := \sqrt{m}+\sqrt{n} - \sqrt{k}$
and $0 < |\t| < \hf$. Then we have
$$
m + n + 2\sqrt{mn} = (\sqrt{m}+\sqrt{n})^2 =
(k + \t)^2 = k + 2\t\sqrt{k} + \t^2.
$$
If $4mn$ is a square, then $2\t\sqrt{k} + \t^2$ is an integer, say $v$.
If $v=0$, then $\t=0$, which is impossible. If $v\not=0$, then
$1 \le |v| \le 2|\t|(\sqrt{k}+1)$, implying $|\t| \gg 1/\sqrt{k}$, which
is better than (2.2). If $4mn$ is not a square, then we have
$$
\t(2\sqrt{k}+\t) = P + 2\sqrt{mn},\quad P = k-m-n \in \ZZ.
$$
Hence, by (2.4),
$$
|\t|\sqrt{k} \;\gg\; ||2\sqrt{mn}|| \;\gg\; {1\over\sqrt{mn}},
$$
giving (2.2).

The proof of (2.3) is on the same lines as the proof of (2.2), only
it is a little more involved.
It suffices to consider only the case of the `+' sign,
since both cases are analogous. Hence let
$$
\rho  := \sqrt{m}+\sqrt{n} - \sqrt{k} - \sqrt{l} \leqno(2.5)
$$
with $\rho \not= 0$. We can
obviously suppose that $|\rho| \le 1/(5\sqrt{k})$. Squaring (2.5) it
follows that $2\sqrt{mn} - 2\sqrt{kl} = v + \rho\mu$, where $v =
k+l-m-n$ is an
integer and $\mu = 2(\sqrt{k} + \sqrt{l}) + \rho$, so that
$0 < \mu \le 4\sqrt{k} + \mu < 5\sqrt{k}$. If $v=0$, then
a better result than (2.3) will follow. If $v\not=0$,
another squaring yields
$$
- 8\sqrt{mnkl} = Q + \rho\nu\qquad(Q\in \ZZ),
$$
with $\nu = 2v\mu + \rho\mu^2.$ If $\nu = 0$, then $\rho\mu = 2|v| \ge 2$,
which contradicts $\rho \le 1/(5\sqrt{k})$. Hence $\nu \not=0$. If
$64mnkl$ is  a square, then again a better result than (2.3) will follow.
If this is not the case,  then $|v| = |k+l-m-n| \le 2k$. Therefore we have,
since $0 < \nu \ll k^{3/2}$,
$$
k^{3/2}|\rho| \gg |\rho\nu| \ge ||8\sqrt{mnkl}|| \,\gg\, {1\over\sqrt{mnkl}},
$$
which gives then (2.3) by an obvious analogue of (2.4).

\medskip
Tsang's proof of (1.7) with $\gamma = 45/23$ depended on the following
lemma.
\medskip
LEMMA (K.-M. Tsang [23]). For any real numbers $\a \not = 0, \b$ and
$0 < \delta < \hf$, we have uniformly
$$
\#\,\Bigl\{\, K < k \le 2K\;:\; ||\b + \a\sqrt{k}|| < \delta
\Bigr\} \ll K\delta + |\a|^{1/3}K^{1/2} + |\a|^{-1/2}K^{3/4}.\leqno(2.6)
$$

 \medskip We shall present now a lemma which, in the relevant range
needed for the proof of Theorem 2, supersedes (2.6). This is

\medskip
LEMMA 3. {\it For real numbers $0 < \delta < \hf$, $\b$ and
$\a \gg 1$}
$$
\#\,\Bigl\{\, K < k \le 2K\;:\; ||\b + \a\sqrt{k}|| < \delta
\Bigr\} \ll_\e K\delta + |\a|^{1/2}K^{1/4+\e}
+ K^{1/2+\e},\leqno(2.7)
$$
{\it where the $\ll$-constant depends only on $\e$}.

\medskip
{\bf Proof of Lemma 3.} We may suppose $\a > 0$.
Denote by ${\Cal N}$ the expression on the
left-hand side of (2.7). We shall use an elementary idea
 contained e.g., in  M.N. Huxley [4, Lemma 3.1.1].
This says that the number of integer points close to the
function ($\b + \a\sqrt{k}$ in this particular case) is essentially
the same as the number of points close to the function that is the
inverse of the original function (with appropriate new $\delta$).
 For our problem,
in his notation $T = \a, L \asymp K, F(x) = \sqrt{x} + \b/\a, x \asymp K$,
$\delta' = \delta\sqrt{K}/\a, G(y) = (y - \b/\a)^2$. Then we obtain
$$
{\Cal N} \ll 1 + \delta' + {\Cal N' } \ll 1 +
 \delta\sqrt{K} + {\Cal N' },\leqno(2.8)
$$
where
$$
{\Cal N' } = \#\,\Bigl\{\, n \asymp \a\sqrt{K} \;:\;
||G(n/T)|| < \delta'\Bigr\}.
$$
Now  we may easily reduce the problem
 to the estimation of exponential sums. For example,
[4, Lemma 5.3.2] gives (setting $Y = \a\sqrt{K}$)
$$
{\Cal N}' \ll Y\delta' + YH^{-1}
+ H^{-1}\sum_{h=1}^H\, \Bigl|\sum_{\nu\asymp Y}{\roman e}
\left(hG\bigl({\nu\over T}\bigr)\right)
\Bigr|,
$$
where ${\roman e}(z) = {\roman e}^{2\pi iz}$, and $H$ is an integer
 satisfying $ H \asymp (\delta')^{-1}$.
To the above exponential sum we apply the $A$-process of van der
Corput (see e.g., [6, Lemma 2.5]). For an integer $R$ with
$1 \le R \le Y$ we have
$$\eqalign{\cr&
{\Cal N}' \ll Y\delta' + YH^{-1} +  YR^{-1/2} \cr&
 + Y^{1/2}\left\{(HR)^{-1}\sum_{h=1}^H\,\sum_{r=1}^R\,
\Bigl|\sum_{\nu\asymp Y}{\roman e}
(2hr\nu\a^{-2})\Bigr|\right\}^{1/2}.\cr}
\leqno(2.9)
$$
The sum over $\nu$ is  a geometric
series, and it follows, on  taking
$R\asymp Y$, and supposing $\delta' < 1$
(for otherwise the final result is obvious), that the expression
in curly brackets in (2.9) is
$$
\eqalign{\cr&
\ll_\e (HR)^{\e-1}\sum_{\mu=1}^{2HR}
\min\left({\Bigl|\Bigl|{\mu\over\a^2}\Bigr|\Bigr|}^{-1},\,Y\right)\cr&
\ll_\e (HR)^{\e-1}(1 + HR\a^{-2})\left(\sum_{1\le \mu\le{1\over2}\a^2}
\min\left({\Bigl|\Bigl|{\mu\over\a^2}\Bigr|\Bigr|}^{-1},
\,Y\right)+Y\right),\cr}
\leqno(2.10)
$$
since $||x|| = ||1-x||$, and $||\mu/\a^2||=0$ if $\mu = \a^2\,(\in \NN)$.
If $\a \ll \sqrt{K}$, then (for $1 \le \mu \le \hf\a^2$)
$$
{\Bigl|\Bigl|{\mu\over\a^2}\Bigr|\Bigr|}^{-1} = {\a^2\over\mu} \le \a^2 \ll Y,
$$
and if $\a \gg \sqrt{K}$, then
$$
\min\left({\Bigl|\Bigl|{\mu\over\a^2}\Bigr|\Bigr|}^{-1},\,Y\right) \ll Y
$$
for $\mu \le \a^2/Y$. Hence in any  case the second sum in (2.10) is
$$
\ll \sum_{\mu\le{1\over2}\a^2}{\a^2\over\mu} + \sum_{\mu\le\a^2/Y}Y
\ll \a^2\log K.
$$
This gives that the expression
in curly brackets in (2.9) is
$$\eqalign{\cr&
\ll_\e K^\e(1 + \a^2(HR)^{-1} + Y\a^{-2} + Y(HR)^{-1})\cr&
\ll_\e K^\e(1 + \a^2\delta'Y^{-1} + Y\a^{-2} +  \delta')\cr&
\ll_\e K^\e(1 + \delta' + Y\a^{-2}).\cr}
$$
Inserting this bound in (2.9), we obtain
$$
{\Cal N}' \ll_\e Y\delta' +K^\e(Y^{1/2}( 1 + (\delta')^{1/2}) + Y/\a)
\ll_\e K\delta + K^{1/4+\e}\a^{1/2} + K^{1/2+\e},
$$
and (2.7) follows then from (2.8) and (2.9).

\medskip
LEMMA 4. {\it Let ${\Cal N}$ denote the number of solutions  in integers
m,n,k of the inequality}
$$
|\sqrt{m} + \sqrt{n} - \sqrt{k}| \le \delta\sqrt{M} \qquad(\delta>0)
\leqno(2.11)
$$
{\it with $M' < n \le 2M', M < m \le 2M, k \in \NN$ and $M' \le M$. Then}
$$
{\Cal N} \ll_\e M^{\e}(M^2M'\delta + (MM')^{1/2}).\leqno(2.12)
$$

\medskip
{\bf Proof of Lemma 4}. We have adopted the
condition $M' < n < 2M'$ instead of the more natural one $N < n \le 2N$,
because we did not want to confound ourselves with the parameter $N$ in (3.4).
The result is trivial if $\delta \ge {1\over4}$. If
$\delta < {1\over4}$, then squaring (2.11) we obtain $k \asymp M$ and also
$(\sqrt{m} + \sqrt{n})^2 = k + O(\delta M)$, implying that
$$
\quad 2\sqrt{mn} = P + O(\delta M),\leqno(2.13)
$$
where $P = k - m - n$ is an integer. If $\delta \gg 1/M$, then (for a
given pair $(k,n)$) it follows on squaring (2.11), that there are
$\ll M\delta $ choices for $m$, hence the bound in (2.12) is trivial.
If $0 < \delta < c/M$ for sufficiently small $c>0$, then from (2.13)
it follows that $||2\sqrt{r}|| < \delta'$ with $r = mn \ll MM'$ and
$\delta' = \delta M$. Since for each $r$ there are at most
$d(r) (\ll_\e r^\e)$ choices for $(m,n)$, the bound (2.12) follows from
Lemma 3 (with $\a = 2, \b = 0$).

\medskip
LEMMA 5. {\it Let ${\Cal N}$ denote the number of solutions  in integers
m,n,k,l of the inequality}
$$
|\sqrt{m} + \sqrt{n} - \sqrt{k} - \sqrt{l}|
\le \delta\sqrt{K}\quad(\delta>0)\leqno(2.14)
$$
{\it with the conditions} $M < m \le 2M$, $M' < n \le 2M'$, $K < k \le 2K$,
$L < l \le 2L$, $1 \ll M,M',L \ll K$. {\it Then}
$$
{\Cal N} \ll_\e KMM'L(\delta + K^{\e-3/2}) + K\min(M,M',L),\leqno(2.15)
$$
{\it and also}
$$
{\Cal N} \ll_\e KMM'L\Bigl(\delta K^2 + (KMM'L)^{-1/2}\Bigr)K^\e.
\leqno(2.16)
$$

\medskip
{\bf Proof of Lemma 5}. Lemma 5 is the crucial lemma in the proof of
Theorem 2. It should be compared to Lemma 1 when $k = 2$, in which
case (2.1) provides a sharper bound. However, the variables in Lemma 1
are supposed to be of the same order of magnitude. This condition is not
fulfilled here, and it accounts for several technical difficulties.
The point of two estimates, (2.15) and (2.16), is that for exceptionally
small $\delta$ the bound (2.16) supersedes (2.15), and this fact will
be used in the proof of Theorem 2.

\smallskip
 We begin the proof of (2.15) with some preliminary observations.
If $\delta \gg 1$,
then (2.15) is trivial.  If $c/K \le \delta \ll 1$ for any constant $c>0$
then from (2.14)  we have
$$
k = (\sqrt{m} + \sqrt{n} - \sqrt{l})^2   +  O(\delta K).
$$
This implies that, for a given triplet $(m,n,k)$,
there are at most $\ll \delta K$
choices for $k$, hence trivially there are no more that $\ll MM'\delta KL$
choices for $(m,n,k,l)$, yielding  the first bound in (2.15).
For $\delta < c/K$
and $c$ small we must
have either $M \gg K$ or $M' \gg K$, for otherwise (2.14) is impossible.
Suppose that the former holds, so that $M \asymp K$.
Thus, as far as the order of magnitude is concerned,
the variables $m$ and $k$
play a symmetric r\^ole. Therefore in what concerns the order
of $M'$ and $L$ we
may assume, without loss of generality, that $M' \ll L$.
When  $\delta <  c/K$,
then if $k = n$ we must have $l = m$, and if $k = m$, then $l = n$,
for otherwise
(2.14) cannot hold. In these cases the number of solutions
is $\ll KM'$, which
is accounted by the last term in (2.15). Henceforth
we assume that $k \not = n$,
$k\not = m$.

\medskip
From (2.14) we obtain by squaring
$$
k + l + n + 2\sqrt{kl} - 2\sqrt{kn} - 2\sqrt{ln} = m + O(\delta K).
$$
This implies, since $\delta K <c$ with small $c>0$, that
$$
||\b + \a\sqrt{k}|| \ll \delta K\quad(\a = 2(\sqrt{l} - \sqrt{n}\,)
\not= 0, \,
\b = -2\sqrt{ln}\,).\leqno(2.17)
$$
Denote by $R(l,n)$ the number of solutions (in $K < k \le 2K$) of (2.17).
We have to consider two cases when we estimate $R(l,n)$.

\smallskip
a) {\it The case when} $\sqrt{l} - \sqrt{n} \gg 1$.

\smallskip
In this case we apply Lemma 3 with $ 1 \ll \alpha
=  \sqrt{l} - \sqrt{n} \ll \sqrt{K}$
to obtain that
$$
\sum_{l,n,\sqrt{l} - \sqrt{n} \gg 1}R(l,n)
\ll_\e M'L(K^2\delta + K^{1/2+\e}).\leqno(2.18)
$$

\smallskip
b) {\it The case when} $0 < |\sqrt{l} - \sqrt{n}| \ll 1$.
\smallskip

In this case we apply [4, Lemma 3.1.2]. In the notation of this lemma we have
to take $\delta K$ as  $\delta$, $L = K, T = \a =
\sqrt{l} - \sqrt{n} \,(\ll 1),
 M = \sqrt{K}$. Therefore
$$
R(l,n) \ll 1 + K^2\delta + K^{1/2}\a + K^{3/2}\delta\a^{-1}.\leqno(2.19)
$$
Set $r = |n-l|$, so that $r>0$, since $\sqrt{l} - \sqrt{n}\not=0$.
If $l > n$, then
$l = n + r$ with $r \ll L$, while if $n > r$, then $n = l + r$ with
$r \ll M'$. Thus using (2.19) we obtain
$$\eqalign{\cr&
\sum\limits_{l,n,0<|\sqrt{l} - \sqrt{n}| \ll 1}R(l,n)
\ll_\e L\sum_{1\le r\ll\sqrt{M'}}
(1 + K^2\delta + K^{1/2}\a + K^{3/2}\delta\sqrt{M'}r^{-1})\cr&
\qquad + M'\sum_{1\le r\ll\sqrt{L}}
(1 + K^2\delta + K^{1/2}\a + K^{3/2}\delta\sqrt{L}r^{-1})\cr&
\qquad\ll K^2L(M')^{1/2}\delta + K^{1/2}L(M')^{1/2}
+ LK^{3/2}(M')^{1/2}\delta\log K\cr&
\qquad+ K^2M'L^{1/2}\delta + K^{1/2}M'L^{1/2}
+ M'K^{3/2}L^{1/2}\delta\log K\cr&
\qquad\ll_\e M'L(K^2\delta + K^{1/2+\e}).
\cr}\leqno(2.20)
$$
From (2.18) and (2.20) we finally obtain, since $M\asymp K$ and $M' = \min(M,M',L)$,
$$
\eqalign{\cr
{\Cal N} &\ll KMM'L\delta + KM' + \sum_{l,n}R(l,n)\cr&
\ll_\e  KMM'L\delta + KM' + M'LK^2\delta + M'LK^{1/2+\e}\cr&
\ll_\e  KMM'L(\delta + K^{\e-3/2}) + K\min(M,M',L).\cr}
$$

\medskip
It remains yet to prove the bound (2.16) of Lemma 5. Proceeding from (2.5),
as in the proof of (2.3), we have that
$$
8\sqrt{mnkl} = Q + O(|\rho|K^{3/2})\qquad(Q\in\ZZ),
$$
where $\rho \,(\ll \sqrt{K})$ is given by (2.5), which may be written as
$$
8\sqrt{j} = Q + O(\delta K^2)\qquad(j = mnkl).
$$
The bound (2.16) follows then from Lemma 3 with $\a = 8$, since for
each given $j$ there are $\ll_\e K^\e$ choices of quadruples $(m,n,k,l)$.
This completes the proof of Lemma 5.

\medskip
LEMMA 6. {\it Let ${\Cal N}$ denote the number of solutions  in integers
m,n,k,l of the inequality}
$$
0 < |\sqrt{m} + \sqrt{n} + \sqrt{k} - \sqrt{l}|
\le \delta\sqrt{K}\quad(\delta>0)
\leqno(2.21)
$$
{\it with the conditions} $M < m \le 2M$, $M' \le N \le 2M'$, $K < k \le 2K$,
$L < l \le 2L$, $1 \ll M,M' \ll K$. {\it Then}
$$
{\Cal N} \ll_\e KMM'L(\delta + K^{\e-3/2}),\leqno(2.22)
$$
{\it and also}
$$
{\Cal N} \ll_\e KMM'L(\delta K^2 + (KMM'L)^{-1/2})K^\e.\leqno(2.23)
$$

\medskip
{\bf Proof of Lemma 6}. The proof is on the same lines
as the proof of Lemma 5,
only it is  less difficult, and the details will be
therefore omitted. It is only
the case of $|\sqrt{m} + \sqrt{n} + \sqrt{k} - \sqrt{l}|$
that is treated in [23],
while the case of the more difficult problem with
$|\sqrt{m} + \sqrt{n} - \sqrt{k} - \sqrt{l}|$ is only
mentioned on top of p. 76: ``By
a similar argument, we show that the same estimate holds
 with $S_4(x)$ in place
of $S_5(x)$.'' Although this is essentially  true, it is
the proof of the more difficult case that should have been given.

To get back to Lemma 6, note that the
condition $1 \ll M,M' \ll K$ may be assumed without loss of generality. However,
then from (2.14) we have
$$
l  = (\sqrt{m} + \sqrt{n} + \sqrt{k})^2 + O(\delta K),
$$
since $\delta < c/K$ may me assumed as in the case of Lemma 5. This gives
$L \asymp K$, and then the proof is analogous to the proof of Lemma 5,
but easier. The term $K\min(M,M',L)$, present in (2.15), is not
necessary in (2.22).  Finally the proof of (2.23)
is completely analogous to the proof of (2.16).

\head
3.  Proof of Theorem 1
\endhead

For the sake of simplicity and readability we shall retain, both in the proof
of Theorem 1 and Theorem 2, the notation of [23].
We shall use a modified form of (2.11)--(2.12), due to T. Meurman [9],
which shows that, for most $x$, the partial sum approximating $\D(x)$
has a small error, provided that the length of the sum is sufficiently
large. More precisely, [9, Lemma 3] says that, for $Q \gg x \gg1$,
$$
\D(x) = {\txt{x^{1/4}\over\pi\sqrt{2}}}\sum_{n\le Q}
d(n)n^{-3/4}\cos(4\pi\sqrt{nx}-
{\txt{1\over4}}\pi) + F(x) = {\txt{1\over\pi\sqrt{2}}}
\sum\nolimits_Q(x) + F(x),\leqno(3.1)
$$
where $F(x) \ll x^{-1/4}$ if $||x|| \gg x^{5/2}Q^{-1/2}$ and otherwise
$F(x) \ll_\e x^\e$. In (3.1) we take $Q = H^7$. Then we have
$$
\int_H^{2H}(\D(x)- F(x))^3\d x = ({\txt{1\over\pi\sqrt{2}}})^3\int_H^{2H}
\sum\nolimits_Q(x)^3\d x. \leqno(3.2)
$$
The left-hand side of (3.2) equals
$$\eqalign{\cr&
\int_H^{2H}(\D(x)+ O(H^{-1/4}))^3\d x
- \int_{H,||x||\ll1/H}^{2H}(\D(x) + O_\e(H^\e))^3\d x\cr&
= \int_H^{2H}\D^3(x)\d x + O\left(H^{-1/4}\int_H^{2H}\D^2(x)\d x\right)
+ O_\e(H^\e\sup_{H\le x\le 2H}|\D(x)|^3)\cr&
= \int_H^{2H}\D^3(x)\d x + O(H^{5/4}),\cr}
$$
since $\D(x) \ll x^{1/3}$ and $\int_H^{2H}\D^2(x)\d x \ll H^{3/2}$
(see [6, Chapter 13]). Now in (3.2) we write
$$
\sum\nolimits_Q(x) = \sum\nolimits_N(x) + R_{N,Q}(x)
\qquad(H \le x \le 2H,\;1\ll N \ll H,\; Q = H^7)
$$
with
$$
\eqalign{\cr
\sum\nolimits_N(x)  &= x^{1/4}\sum_{n\le N}d(n)n^{-3/4}\cos(4\pi\sqrt{nx}
- {\txt{1\over4}}\pi),\cr
R_{N,y}(x) & = x^{1/4}\sum_{ N< n \le y}d(n)n^{-3/4}\cos(4\pi\sqrt{nx}
- {\txt{1\over4}}\pi) \cr&
\ll_\e x^\e\left(1+ x^{1/2}N^{-1/2}\right)\quad(N < y \ll x^A).
\cr}\leqno(3.3)
$$
The bound for $R_{N,y}(x)$ in (3.3) follows from (1.14), when we write
it once with $N$ and once with $y$ replacing $N$, and then subtract
the resulting expressions. It follows that  (3.1) gives
$$\eqalign{\cr&
\int_H^{2H}\D^3(x)\d x = {1\over(\pi\sqrt{2})^3}\int_H^{2H}
\sum\nolimits_N^3(x)\d x + O(H^{5/4})
\cr& +
{1\over(\pi\sqrt{2})^3}\int_H^{2H}\left(3\sum\nolimits_N^2(x)R_{N,Q}(x) +
 3\sum\nolimits_N(x)R_{N,Q}^2(x) + R_{N,Q}^3(x)\right)\d x.\cr}\leqno(3.4)
$$
Had we used directly (1.14), we would have obtained (3.4), but
with the error term $O_\e(H^{3/2+\e})$ instead of $O(H^{5/4})$, and
this is too large for our purposes.

\smallskip
We shall first evaluate the integral with $\sum_N^3(x)$,
and then the remaining
integrals. If we used directly the bound (3.3) for $R_{N,y}(x)$,
 we would obtain the additional error terms (cf. [23])
$O_\e(H^{2+\e}N^{-{1\over2}})
+ O_\e(H^{{5\over2}+\e}N^{-{3\over2}})$, which would be too large
to yield the exponent $\b = 7/5$ in the error term in Theorem 1.
Put
$$
r = r(m,n,k) := d(m)d(n)d(k)(mnk)^{-3/4}\qquad( 1 \le m,n,k \le N),
\leqno(3.5)
$$
and $r=0$ otherwise. Then
$$\eqalign{\cr
{\sum\nolimits_N(x)}^3 &= {\txt{3\over4}}\sum rx^{3/4}
\cos(4\pi(\sqrt{m}+\sqrt{n}-\sqrt{k})\sqrt{x} - {\txt{1\over4}}\pi)\cr&
+ {\txt{1\over4}}\sum rx^{3/4}
\cos(4\pi(\sqrt{m}+\sqrt{n}+\sqrt{k})\sqrt{x}-{\txt{3\over4}}\pi)\cr&
= S_0(x) + S_1(x) + S_2(x),\cr}
$$
say, where
$$\eqalign{\cr
 S_0(x) &:= {\txt{3\over4\sqrt{2}}}\sum_{\sqrt{m}+\sqrt{n}=
\sqrt{k}}rx^{3/4},\cr
  S_1(x) &:=  {\txt{3\over4}}\sum_{\sqrt{m}+\sqrt{n}\not=\sqrt{k}}rx^{3/4}
\cos(4\pi(\sqrt{m}+\sqrt{n}-\sqrt{k})\sqrt{x}-{\txt{1\over4}}\pi)\cr
 S_2(x) &:=  {\txt{1\over4}}\sum rx^{3/4}
\cos(4\pi(\sqrt{m}+\sqrt{n}+\sqrt{k})\sqrt{x}-{\txt{3\over4}}\pi).\cr}
$$
Tsang [23] has shown that, with explicit $c_1>0$,
$$
\int_H^{2H}S_0(x)\d x = {3c_1\over4\sqrt{2}}\int_H^{2H}x^{3/4}\d x
+ O_\e(H^{{7\over4}+\e}N^{-1})\leqno(3.6)
$$
and, by using the first derivative test, it follows that
$$
\int_H^{2H}S_2(x)\d x \ll_\e H^{{5\over4}+\e}N^{{1\over4}}.\leqno(3.7)
$$
The most delicate task is the estimation of the integral of $S_1(x)$, where
we shall proceed  differently than in [23]. In integrating $S_1(x)$
we have that $E \not=0$, where
$$
E \;=\; E(m,n,k) \;:=\; \sqrt{m}+\sqrt{n}-\sqrt{k}.
$$
But if $E \not=0$, then by (2.2) (Lemma 2) we have that
$$
|E| \gg (mnk)^{-1/2}.\leqno(3.8)
$$
Without loss of generality we may assume that $n\le m$, hence it suffices
to bound
$$
{\sum}^{*}r\int_H^{2H}x^{3/4}\cos(4\pi E\sqrt{x}-{\txt{3\over4}}\pi)
\d x,\leqno(3.9)
$$
where ${\sum}^{*}$ denotes summation over $m,n,k \le N$
such that $n\le m$ and
(3.8) holds. We further suppose that $M' < n\le2M'$,
$M < m \le 2M$, $K < k\le 2K$, so that $M' \ll M$. Then either
$K \ll M$ or $M \ll K$. Suppose first $K \ll M$.
We distinguish the following cases when we estimate (3.9).

\medskip
a) {\it The case} $|E| \gg  m^{{1\over2}}$. Then   the  integral in (3.9) is
estimated, by the first derivative test, as $\ll H^{5/4}|E|^{-1}$.
The corresponding portion of the sum in (3.9) is
$$\eqalign{\cr&
\ll_\e H^{{5\over4}+\e}\sum_{M' < n\le2M'}\sum_{M < m\le2M}\sum_{K < k\le2K}
(mnk)^{-3/4}m^{-1/2}\cr&
\ll_\e H^{{5\over4}+\e}(M'MK)^{1/4}M^{-1/2}\cr&
\ll_\e H^{{5\over4}+\e}(M')^{1\over4}\cr&
\ll_\e H^{{5\over4}+\e}N^{1\over4}.\cr}
\leqno(3.10)
$$

\medskip
b) {\it The case} $ |E| \le  cm^{{1\over2}}$ with small $c>0$.
Then it follows that
$$
k = (\sqrt{m}+\sqrt{n})^2 + O(|E|\sqrt{m})  = (\sqrt{m}+\sqrt{n})^2 + O(cm),
$$
hence $k \asymp m$ if $c$ is sufficiently small.
 We use the first derivative test and Lemma 4 (with
$\delta = |E|\sqrt{M}$) to obtain that the contribution will be
$$
\eqalign{\cr&
\ll_\e H^{{5\over4}+\e}\max_{M'\le M\ll N, |E|\gg 1/\sqrt{M^2M'}}
(M^2M')^{-{3\over4}}|E|^{-1}(M^{3\over2}M'|E| + (MM')^{1\over2})\cr&
\ll_\e H^{{5\over4}+\e}N^{1\over4} +  H^{{5\over4}+\e}\max_{M'\le M\ll N}
(M^2M')^{-3/4}(M^2M')^{1/2}(MM')^{1/2})\cr&
\ll_\e H^{{5\over4}+\e} N^{1\over4}.\cr}
$$
If $M \ll K$, we repeat the above argument according to the cases when
$|E| \gg  k^{{1\over2}}$ and $ |E| \le  ck^{{1\over2}}$. The bound
corresponding to (3.10) will be, in this case,
$$\eqalign{\cr&
\ll_\e H^{{5\over4}+\e}(M'MK)^{1/4}K^{-1/2}
\cr&\ll_\e H^{{5\over4}+\e} (M')^{1\over4}
\ll_\e H^{{5\over4}+\e} N^{1\over4}.\cr}
$$

\medskip
It remains to estimate  the integrals in (3.4) with $R_{N,Q}(x)$,
all of which will contribute to the error terms. Namely when
we expand each of the three terms, we shall obtain sums analogous
to $S_0(x), S_1(x), S_2(x)$, but at least one variable, say $k$, will
satisfy $k > N$. Should it happen that a sum of two roots equals
a third root, then this can happen only if $\sqrt{m} + \sqrt{n}
= \sqrt{k}$. But then the argument of [23], which yielded (3.6),
clearly shows that the contribution of such sums must be
$\ll_\e H^{7/4+\e}N^{-1}$. In the case of the integrals of
$\sum_N(x)R_{N,Q}^2(x)$ and $\sum^2_N(x)R_{N,Q}(x)$,
the analogues of the sums $S_1(x)$
and $S_2(x)$ are estimated as $\ll_\e H^{{5\over4}+\e} N^{1\over4}$
by the above method of proof, since we have $M' \ll N$, and this is
the crucial condition both in case a) and case b).
However, when we deal with the integral of $R^3_{N,Q}(x)$, this approach
does not work, since all variables are $\gg N$. Instead we use
$R_{N,Q}(x) \ll_\e H^{1/2+\e}N^{-1/2}$ and the approach used in the
estimation of the integral $I_2$ in (2.13). This yields,
on using the first derivative test,
$$\eqalign{\cr&
\int_H^{2H}R_{N,Q}^3(x)\d x
\ll_\e H^{1/2+\e}N^{-1/2}\int_H^{2H}R_{N,Q}^2(x)\d x\cr&
\ll_\e H^{1/2+\e}N^{-1/2}\int_H^{2H}H^{1/2}\Bigl|
\sum_{N<n\le H^7}d(n)n^{-3/4}
\cos(4\pi\sqrt{nx}-{\txt{1\over4}}\pi)\Bigr|^2\d x\cr&
\ll_\e H^{1/2+\e}N^{-1/2}\left(H^{3/2}\sum_{n>N}d^2(n)n^{-3/2} +
H\sum_{N<m\not=n<H^7}{d(m)d(n)\over(mn)^{3/4}|\sqrt{m}-\sqrt{n}|}
\right)\cr&
\ll_\e H^{1/2+\e}N^{-1/2}(H^{3/2}N^{-1/2} + H)\cr&
\ll_\e H^{2+\e}N^{-1},\cr}\leqno(3.11)
$$
since $N \ll H$. Thus gathering all
the estimates we arrive at
$$
\int_H^{2H}\D^3(x)\d x =
{3c_1\over4\sqrt{2}}\int_H^{2H}x^{3/4}\d x + O_\e(H^{2+\e}N^{-1})
+ O_\e(H^{{5\over4}+\e}N^{1\over4}).\leqno(3.12)
$$
The proof is completed when we take $N = H^{3/5}$, then $H = X/2, X/2^2,
\ldots\,$ and sum all the resulting expressions. Note that
it is only in (3.11) that the error term $H^{2+\e}N^{-1}$ appears,
while before we had $H^{7/4+\e}N^{-1}$. If we could obtain (3.12)
with the latter error term instead of  $H^{2+\e}N^{-1}$, we
would then obtain  the exponent $\b = 27/20$ in (1.6). This would be
then  the limit of the present method.

\head
4.  Proof of Theorem 2
\endhead
We pass now to the proof of Theorem 2. As in the proof of Theorem 1,
our approach differs from Tsang's in the treatment of the critical
sums $S_4, S_5$ below, and additional saving
comes since we use Lemma 5 and Lemma 6
instead of Tsang's (2.7). Although slight improvements
of his result could be obtained by using the theory of exponent pairs to
estimate the exponential sum appearing on the right-hand side of
[23, eq. (4.9)], these improvements cannot attain the strength of our
Lemma 5 and Lemma 6.

\smallskip
 We start from (1.14) to obtain
$$
\eqalign{\cr&
\int_{H}^{2H}\D^4(x)\d x = (\pi\sqrt{2})^{-4}
\int_{H}^{2H}{\sum\nolimits_N(x)}^4\d x \cr&
+ O_\e\left(H^{{1\over2}+\e}N^{-{1\over2}}\int_{H}^{2H}
{|\sum\nolimits_N(x)|}^3\d x\right)
+ O_\e(H^{3+\e}N^{-2})\cr&
= (\pi\sqrt{2})^{-4}
\int_{H}^{2H}{\sum\nolimits_N(x)}^4\d x +  O_\e(H^{{9\over4}+\e}N^{-1/2})
+ O_\e(H^{3+\e}N^{-2}),\cr}
\leqno(4.1)
$$
with $\sum_N(x)$ given by (3.3). Here, unlike in (3.4), we used
the crude estimate $R_{N,H}(x) \ll_\e x^{1/2+\e}N^{-1/2}$ and the bound
$\int_{H}^{2H}|\D^3(x)|\d x \ll_\e H^{7/4+\e}$
(see [5] or [6]). The error terms in (4.1) suffice for (1.7)
with $\gamma = 23/12$, although they could be improved by a
technique similar to the one used in the proof of Theorem 1.

\smallskip
To evaluate the integral with $\sum_N^4(x)$ in (4.1),
 we proceed similarly as in Tsang [23]. With
$$
r_1 = r_1(m,n,k,l) := (mnkl)^{-3/4}d(m)d(n)d(k)d(l)
\quad(m,n,k,l \le N; m,n,k,l \in \NN)
$$
and $r_1 \equiv 0$ otherwise, we have
$$
\eqalign{\cr
{\sum\nolimits_N(x)}^4 &= S_3(x) + S_4(x) + S_5(x) + S_6(x),\cr
S_3(x) &= {\txt{3\over2}}\sum_{\sqrt{m} + \sqrt{n} =
\sqrt{k}+\sqrt{l}}r_1x,\cr
S_4(x) &= {\txt{3\over8}}\sum_{\sqrt{m} + \sqrt{n} \not =
\sqrt{k}+\sqrt{l}}r_1x
\cos(4\pi(\sqrt{m} + \sqrt{n} - \sqrt{k}-\sqrt{l})\sqrt{x}\,),\cr
S_5(x) &= {\txt{1\over2}}\sum r_1x
\sin(4\pi(\sqrt{m} + \sqrt{n} + \sqrt{k}-\sqrt{l})\sqrt{x}\,),\cr
S_6(x) &= -{\txt{1\over8}}\sum r_1x
\cos(4\pi(\sqrt{m} + \sqrt{n} + \sqrt{k}+\sqrt{l})\sqrt{x}\,).\cr}
$$
As in [23] it follows that, with suitable $c_2>0$, we have
$$
\eqalign{\cr
\int_{H}^{2H}S_3(x)\d x &= {\txt{3\over8}}c_2
\int_{H}^{2H}x\d x + O_\e(H^2N^{-1/4+\e}).\cr}
$$
By using the first derivative test  we obtain
$$
\int_{H}^{2H}S_6(x)\d x \ll_\e H^{3/2+\e}N^{1/2}.
$$
It remains to consider the integrals with $S_4(x)$ and $S_5(x)$.
To do this, set
$$
\D_\pm = \D_\pm(m,n,k,l) := 4\pi(\sqrt{m} + \sqrt{n} \pm  \sqrt{k}-\sqrt{l}).
$$
We may assume that ($M \ll K, M' \ll K$)
$$
M < m \le 2M,\; M' < n \le 2M', \;K < k \le K' \le 2K,\; L < l \le L'\le 2L
\leqno(4.2)
$$
 and that
$\D_\pm \not = 0$. In the case of $\D_-$ the condition $\D_- \not=0$ is
assumed in $S_4(x)$, and if $\D_+ = 0$, then $S_5(x)$ vanishes identically.
If $\D_\pm \not = 0$, then by (2.3) we have
$$
 |\D_\pm|\; \gg \;k^{-2}(mnl)^{-1/2}\qquad(m \le k,\, n\le k).\leqno(4.3)
$$
We suppose that
$$
|\D_\pm|\; \asymp \; \delta\sqrt{K},\leqno(4.4)
$$
and (similarly as in the proof of Lemma 5)
we may assume that $0 < \delta \ll 1/K$. Namely,
if $\delta \gg 1$, then the number of solutions of (4.4)
is trivially $\ll MM'KL\delta$, hence by the first derivative test and
trivial estimation the relevant portions of $S_4(x)$ and $S_5(x)$ are bounded by
$$
H^{{3\over2}+\e}|\D_\pm|^{-1}(MM'KL)^{-3/4}MM'KL\delta
\ll_\e H^{{3\over2}+\e}N^{1\over2},
\leqno(4.5)
$$
since $ K\ll N$. If $1/K \ll
\delta \ll 1$, then $|\D_\pm| \gg K^{-1/2}$, and moreover we have
$$
(\sqrt{m} + \sqrt{n} - \sqrt{l})^2 = (\pm\sqrt{k})^2 + O(\delta \sqrt{kK})
= k +  O(\delta K),
$$
implying that, for a given triplet $(m,n,l)$,
there are at most $\ll \delta K$
choices for $k$. Hence trivially there are no more that $\ll MM'\delta KL$
choices for $(m,n,k,l)$, and  this gives us again the bound in (4.5).
Then we obtain, squaring the defining relation of $\D_\pm$, that
$M \asymp K$.

\smallskip
To bound the integrals of $S_4(x)$ and $S_5(x)$ it will be
sufficient to estimate
$$
\int_{H}^{2H}\mathop{\sum\nolimits^*}\limits_{m,n,k,l}r_1x\,
{\roman e}(\D_\pm\sqrt{x}\,) \d x,\leqno(4.6)
$$
where  ${\sum}^*$ means that (4.2), (4.3) and  the conditions of Lemma 5
(or Lemma 6) hold with $0 < \delta \ll 1/K$.

\smallskip
We shall consider the case $\D = \D_-$ in detail, since the case of
$\D_+$ is analogous, but less difficult, due to the absence of the
last term in (2.15) of Lemma 5 in the bound (2.22) of Lemma 6. We shall
estimate the integral in (4.6) trivially, or by the first derivative test
to obtain
$$
\int_H^{2H}x{\roman e}(\D\sqrt{x}\,)\d x \;\ll\;H^2\min
\left(1,\,{1\over|\D|\sqrt{H}}\right).\leqno(4.7)
$$
After that, we estimate the remaining portion of the sum by Lemma 5.
Depending on which term in (2.15) (or (2.16)) dominates, we may use
either of the bounds in (4.7) as we please. We consider two cases.

\smallskip
a) {\it The case when} $K \gg H^{1/6}$. Then by using (2.15)
(with $\delta \asymp |\D|K^{-1/2}$) and (4.7)
we see that the relevant portion of (4.6) is
$$
\eqalign{\cr&
\ll_\e H^{2+\e}\max_{M\asymp K,M'\ll K,H^{1/6}\ll K\ll N}
{\sum_{m,n,k,l}}^{*} r_1\min
\Bigl(1,\,{1\over|\D|\sqrt{H}}\Bigr)\cr&
\ll_\e H^{2+\e}\max_{M\asymp K, K\gg H^{1/6}}
\bigl\{(KMM'L)^{1/4}((KH)^{-1/2} + K^{-3/2})\cr&
\quad + (KMM'L)^{-3/4}K(M'L)^{1/2}\bigr\}\cr&
\ll_\e H^{2+\e}N^{1/2} + \max_{K\gg H^{1/6}}H^{2+\e}K^{-1/2}
\cr&
\ll_\e H^{2+\e}N^{1/2} + H^{23/12+\e}.\cr}\leqno(4.8)
$$

\smallskip
b) {\it The case when} $K \ll H^{1/6}$. Then by using (2.16)
(with $\delta \asymp |\D|K^{-1/2}$) and (4.7)
we see that the relevant portion of (4.6) is
$$
\eqalign{\cr&
\ll_\e H^{2+\e}\max_{M\asymp K,M'\ll K,K\ll H^{1/6}}
{\sum_{m,n,k,l}}^{*} r_1\min
\Bigl(1,\,{1\over|\D|\sqrt{H}}\Bigr)\cr&
\ll_\e H^{2+\e}\max_{M\asymp K,K\ll H^{1/6},{\roman (4.3)}}
\Bigl\{(MM'KL)^{1\over4}(K^{-{1\over2}}K^2H^{-1\over2} +
{(MM'KL)^{-{1\over2}}\over|\D|\sqrt{H}}\Bigr\}\cr&
\ll_\e H^{2+\e}\max_{K\ll H^{1/6},{\roman (4.3)}}\Bigl(H^{-1/2}K^{5/2} +
H^{-1/2}|\D|^{-1}(KMM'L)^{-1/4}\Bigr)
\cr&
\ll_\e H^{2+\e}\Bigl(H^{-1/12} +
\max_{K\ll H^{1/6},{\roman (4.3)}}
H^{-1/2}|\D|^{-1}K^{-1/2}(M'L)^{-1/4}\Bigr) .\cr}\leqno(4.9)
$$
But as $\D \not=0$, (4.2) and (4.3) give
$$
|\D|^{-1} \;\ll\; K^{5/2}(M'L)^{1/2}
$$
and therefore the contribution of (4.9) will be
$$
\ll_\e H^{2+\e}\Bigl(H^{-1/12} + \max_{K\ll H^{1/6}} K^{5/2}H^{-1/2}
\Bigr) \ll_\e H^{23/12+\e}.
$$

Therefore putting together all the estimates, we obtain
$$
\eqalign{\cr &
\int_{H}^{2H}\D^4(x)\d x =
{\txt{3\over8}}c_2\int_{H}^{2H}x\d x  + O_\e(H^{2+\e}N^{-{1\over4}}) +
 O_\e(H^{{3\over2}+\e}N^{1\over2})\cr&
+  O_\e(H^{{9\over4}+\e}N^{-{1\over2}}) +
O_\e(H^{3+\e}N^{-2}) +  O_\e(H^{23/12+\e}).
\cr}\leqno(4.10)
$$
The choice $N = H^{3/4}$ yields the assertion of Theorem 2, when
we take $H = X/2, X/2^2, \ldots\,$ and sum the
resulting expressions. The limit of the method is the exponent
$\gamma = 11/6$ in (1.7). This would follow, with $N = H^{2/3}$, if all
the error terms in (4.10) could be made to be $ O_\e(H^{2+\e}N^{-{1\over4}}) +
 O_\e(H^{{3\over2}+\e}N^{1\over2})$.

\head
5.  Proof of Theorem 3
\endhead
The starting point for both (1.10) and (1.11) is the formula
$$
\pi\sqrt{2}\D(x) = \sum\nolimits_N(x) + R_{N,X}(x) + O_\e(x^\e)
\quad(X\le x\le X+H, N \ll X),\leqno(5.1)
$$
which follows on combining (1.14) and (3.3), and where the parameter $N$
will be suitably chosen. There are many details in the proof similar
to the proof of Theorem 1 and Theorem 2, so that we may be fairly brief.
From (5.1) we obtain
$$
\int_X^{X+H}\D^3(x)\d x = (\pi\sqrt{2})^{-3}
\int_X^{X+H}\sum\nolimits^3_N(x)\d x + O({\Cal R}_1+{\Cal R}_2+
{\Cal R}_3+HX^{2/3}),\leqno(5.2)
$$
say, where we used (5.6) with $\t < 1/3$ and we set
$$
{\Cal R}_1 := \Bigl|\int_X^{X+H}\sum\nolimits^2_N(x)R_{N,X}(x)\d x\Bigr|,
$$
$$
{\Cal R}_2 := \int_X^{X+H}|\sum\nolimits_N(x)|R^2_{N,X}(x)\d x,\,
{\Cal R}_3 := \int_X^{X+H}|R^3_{N,X}(x)|\d x.
$$
Analogously as in the proof of Theorem 1 we obtain
$$\eqalign{
\int_X^{X+H}\sum\nolimits^3_N(x)\d x &= B\left((X+H)^{7/4} - X^{7/4}\right)
\cr&
+ O_\e(X^{5/4+\e}N^{1/4}) + O_\e(HX^{3/4+\e}N^{-1}).\cr}\leqno(5.3)
$$
Note that trivial estimation  gives
$$
\sum\nolimits_N(x) \ll_\e X^{1/4+\e}N^{1/4}\qquad(X \le x \le X+H),
\leqno(5.4)
$$
hence from (5.1) we obtain
$$
R_{N,X}(x) \ll_\e X^{\t+\e} + X^{1/4+\e}N^{1/4}\qquad(X \le x \le X+H),
\leqno(5.5)
$$
where $\t$ is such a constant ($1/4 \le \t \le 1/3$) for which one has
$$
\D(x) \;\ll_\e\;x^{\t+\e}.\leqno(5.6)
$$
We have, by using the estimation of (3.11) and (5.4),
$$\eqalign{
{\Cal R_2} &\ll \max_{X\le x\le X+H}\Bigl|\sum\nolimits_N(x)\Bigr|
\int_X^{X+H}R_{N,X}^2(x)\d x\cr&
\ll_\e X^{5/4+\e}N^{1/4} + HX^{3/4+\e}N^{-1/4}.\cr}\leqno(5.7)
$$
In a similar way, by using (5.5), it follows that
$$\eqalign{
{\Cal R_3} &\ll \max_{X\le x\le X+H}\Bigl|R_{N,X}(x)\Bigr|
\int_X^{X+H}R_{N,X}^2(x)\d x\cr&
\ll_\e X^\e(X^{5/4}N^{1/4} + HX^{3/4}N^{-1/4} + X^{1+\t}
+ HX^{1/2+\t}N^{-1/2}).\cr}\leqno(5.8)
$$
Developing the integrand in the expression for ${\Cal R}_1$ and using
the first derivative test, it is seen that
$$
{\Cal R}_1 \ll_\e X^{5/4+\e}\sum_{n_1,n_2\le N}\sum_{N<m\le X}
(n_1n_2m)^{-3/4}(\sqrt{m} + \sqrt{n_1} - \sqrt{n_2}\,)^{-1}.
$$
Set $\s := \sqrt{m} + \sqrt{n_1} - \sqrt{n_2}\; ( >0,\,$since
$m > N$ and $  n_1,n_2 \le N)$. The contribution of triplets $(m,n_1,n_2)$
for which $\s \gg \sqrt{m}$ is $\ll_\e X^{5/4+\e}N^{1/4}$. Suppose
now that $\s \asymp \eta\sqrt{m}$, with $\eta >0$ sufficiently small. This
is possible only if $n_2 \asymp N$ and $m \asymp N$. For a fixed $n_1$,
set $m = N + h, n_2 = N - k$ with $h,k\in \NN$. Then $\s \ll \eta\sqrt{m}$
implies
$$
\sqrt{N+h} - \sqrt{N-k} \ll \eta\sqrt{N},
$$
hence $h + k \ll \eta N$. Thus there are $\ll \eta N^2$ choices for
$(m,n_2)$, and the contribution of such $\s$ is again
$\ll_\e X^{5/4+\e}N^{1/4}$.
Replacing $\eta$ by $2^{-j}\eta\;(j\in\NN)$ and noting that there are
$O(\log X)$ values of $j$ by (2.2) of Lemma 2, it follows that
$$
{\Cal R}_1 \ll_\e X^{5/4+\e}N^{1/4}.\leqno(5.9)
$$
If we now choose (recall that $\t \ge 1/4$ must hold; see e.g., [6, Chapter
[13])
$$
N \;=\; X^{4\t-1+\k},\leqno(5.10)
$$
then from (5.2), (5.3) and (5.7)--(5.10) we obtain
$$
\int_X^{X+H}\D^3(x)\d x = B\left((X+H)^{7/4} - X^{7/4}\right) + {\Cal R},
\leqno(5.11)
$$
with
$$
{\Cal R} \ll_\e X^{1+\t+\e} + HX^{3/4+\e-(\t-1/4+\k/4)}.\leqno(5.12)
$$
The classical (trivial) value $\t = 1/3$ gives rise to the exponent 7/12
in (1.10), while (5.11) and (5.12) show that actually 7/12
can be replaced by $1/4+\t$. Thus the conjectural $\t = 1/4$
would replace 7/12 by 1/2 in (1.10), and the value $\t \le 23/73$
(see M.N. Huxley [4]) yields the exponent $165/292 = 0.56506\ldots
< 7/12 = 0.58333\dots\,$.

\medskip
We pass now to the proof of (1.11). We need the following

\bigskip
LEMMA 7. {\it For any given $0 < \k < \hf$ and $X^{2\k} \le N \ll X$
we have}
$$
\int_X^{X+H}R_{N,X}^4(x)\d x \ll_\k X^{5/3+\k} + HX^{1-\k}.\leqno(5.13)
$$

\bigskip
{\bf Proof of Lemma 7}. With $M < M' \le 2M$ we have
$$
\eqalign{&
\int_X^{X+H}R_{N,X}^4(x)\d x \ll \log X\max_{N\le M\ll X}
\int_X^{X+H}R_{M,M'}^4(x)\d x\cr&
\ll_\e X^{1+\e}\max_{N\le M\ll X}\int_{X-H}^{X+H}\f(x)\Bigl|\sum_{M<n\le M'}
d(n)n^{-3/4}{\roman e}^{4\pi i\sqrt{nx}}\Bigr|^4\d x,\cr}\leqno(5.14)
$$
where $\f(x)$ is a smooth, non-negative function supported in $\,[X-2H,
X+2H]\,$ such that $\f(x) = 1$ when $x\in [X-H,X+H]$, and $\f^{(r)}(x)
\ll_r H^{-r}$ for $r = 0,1,2,\ldots\,$. The fourth power of
the above sum equals
$$
\sum_{M<n_1,n_2,n_3,n_4<M'}d(n_1)d(n_2)d(n_3)d(n_4)(n_1n_2n_3n_4)^{-3/4}
{\roman e}^{Di\sqrt{x}},
$$
where
$$
{\Cal D} = {\Cal D}(n_1,n_2,n_3,n_4) := 4\pi(\sqrt{n_1} +
\sqrt{n_2} - \sqrt{n_3} - \sqrt{n_4}\,).
$$
We perform then, in the last integral in (5.14), a large number of
integrations by parts, taking into account that $\f^{(r)}(x) \ll_r H^{-r}$.
It transpires that only the contribution of quadruples $(n_1,n_2,n_3,n_4)$
for which $|{\Cal D}| \le X^{1/2+\e}H^{-1}$ will be non-negligible.
This contribution is estimated, by Lemma 1 (with $k = 2, \delta \asymp
X^{1/2+\e}H^{-1}M^{-1/2}$) and trivial estimation, as
$$
\eqalign{&
\int_{X-H}^{X+H}\f(x)\Bigl|\sum_{M<n\le M'}
d(n)n^{-3/4}{\roman e}^{4\pi i\sqrt{nx}}\Bigr|^4\d x\cr&
\ll_\e HM^{-3}(M^4X^{1/2+\e}H^{-1}M^{-1/2} + M^2)\cr&
\ll_\e HX^{1/2+\e}M^{1/2} + HM^{-1}.\cr}\leqno(5.15)
$$
On the other hand, by (3.3),
$$
R_{M,M'}(x) \ll_\e x^{1/2+\e}M^{-1/2} \qquad(N \le M \ll X).\leqno(5.16)
$$
This gives
$$\eqalign{&
\int_{X}^{X+H}R_{M,M'}^4(x)\d x \ll_\e X^{1+\e}M^{-1}
\int_{X}^{X+H}R_{M,M'}^2(x)\d x\cr&
\ll_\e X^{2+\e}M^{-1} + HX^{3/2+\e}M^{-3/2},\cr}\leqno(5.17)
$$
by the argument used in deriving (3.11).

\medskip
If $N \le M \le X^{1/3+\k}$, then (5.14)-(5.15) give
$$
\int_{X}^{X+H}R_{M,M'}^4(x)\d x \ll_\e X^{5/3+\e+\k/2} + HX^{1+\e-2\k},
$$
while for $M > X^{1/3+\k}$ and $N \ge X^{2\k}$ we infer from (5.17) that
$$
\int_{X}^{X+H}R_{M,M'}^4(x)\d x \ll_\e X^{5/3+\e-\k} + HX^{1+\e-3\k/2}.
$$
The bound in (5.13) follows on combining the last two bounds in conjunction
with (5.14), and taking e.g. $\e = \k/4$ sufficiently small.

\medskip
It is now not difficult to obtain (1.11). From (5.1) we have
$$
\int_X^{X+H}\D^4(x)\d x = (\pi\sqrt{2})^{-4}\int_X^{X+H}
\sum\nolimits_N^4(x)\d x
+ O(\sum_{j=1}^4I_j) + O_\e(HX^{3\t+\e}),\leqno(5.18)
$$
with $\t$ given by (5.6) (for our purposes any $\t < 1/3$ clearly suffices)
and
$$
I_j  := \int_X^{X+H}|R_{N,X}(x)|^j|\sum\nolimits_N(x)|^{4-j}\d x.\leqno(5.19)
$$
By the arguments of Tsang [3] (see also our discussion after (4.1)) we have
$$
\int_X^{X+H}\sum\nolimits_N^4(x)\d x
= (\pi\sqrt{2})^{4}C((X+H)^2 - X^2) + O_\e(HX^{1+\e}N^{-1/4} +
X^{3/2+\e}N^{9/2}),\leqno(5.20)
$$
and in particular, with $N = X^\k$ and sufficiently small $\k >0$,
$$
\int_X^{X+H}\sum\nolimits_N^4(x)\d x \ll HX \qquad(X^{1/2+\delta}
\le H \ll X).\leqno(5.21)
$$
Using Lemma 7, (5.21) and H\"older's inequality for integrals we
obtain, for $j= 1,2,3,4$,
$$
\eqalign{&
I_j \le \left(\int_X^{X+H}R_{N,X}^4(x)\d x\right)^{j/4}
\left(\int_X^{X+H})\sum\nolimits_N^4(x)\d x\right)^{(4-j)/4}\cr&
\ll (X^{5/3+\k/2} + HX^{1-\k/2})^{j/4}(HX)^{1-j/4}\cr&
\ll H^{1-j/4}X^{1+j/6+\k j/8} + HX^{1-\k j/8}.\cr}\leqno(5.22)
$$
If $H \ge X^{2/3+\delta}$ with a fixed $\delta > 0$, then taking
$0 < \k < \delta$ we obtain from (5.22)
$$
I_j \;\ll\; HX^{1-\k/8}\qquad(j = 1,2,3,4).\leqno(5.23)
$$
Inserting (5.20) and (5.23) in (5.18)
we obtain (1.11), with $\k/8$ replacing $\k$.

\medskip
There are possibilities to extend the range of $H$ for which (1.11) holds.
For example, instead of (5.16) we may write explicitly
$$
R_{M,M'}(x) \ll x^{1/4}\Bigl|\sum_{M<n\le M'}d(n)n^{-3/4}
{\roman e}^{4\pi i\sqrt{nx}}\Bigr|,\leqno(5.24)
$$
taking account that the sum in (5.24) reduces to a double exponential
sum on writing $d(n) = \sum_{m_1m_2=n}1$. Estimating the terms $m_1
= m_2$ trivially, we conclude that
$$
\eqalign{
R_{M,M'}(x) &\ll x^{1/4}M^{-1/4}\sum_{1\le m_2\le M'}
\Bigl|\sum_{m_1\in I(m_2)}{\roman e}^{4\pi i\sqrt{m_1m_2x}}\Bigr|\cr&
+ M^{1/8}x^{1/8},\cr}
$$
where $m_1$ runs over an interval $I(m_2)$, which
is contained in an interval of length $\ll M/m_2$.
Since $M$ is of the order of magnitude close to $X^{1/3}$, $m_1
\gg X^{1/6}$, and the fifth derivative of $\sqrt{m_1m_2x}$ with
respect to $m_1$ is sufficiently small, it means
that the fifth derivative test
(see Graham-Kolesnik [1, Th. 2.8] with $q=3$) can be applied. This
application will produce a non-trivial estimation for $R_{M,M'}(x)$ which,
in the range relevant for our problem, will lead to a better value than
2/3 in (1.11).

\vfill
\eject
\topglue2cm
\vskip1cm
\Refs
\bigskip
\item{\bf1.}   S.W. Graham and G. Kolesnilk, `Van der Corput's method
for exponential sums', LMS Lecture Notes Series {\bf 126}, Cambridge
University      Press, Cambridge, 1991.

\item{\bf2.} G.H. Hardy, `On Dirichlet's divisor problem', Proc. London
Math. Soc. {\bf15}(2)(1916), 1-25.

\item{\bf3.} D.R. Heath-Brown, `The distribution of moments
in the Dirichlet divisor problems', {\it Acta Arith.} {\bf60}(1992),
389-415.

\item{\bf4.} M.N. Huxley, `Area, Lattice Points and Exponential
Sums', Oxford Science Publications, Clarendon Press,
Oxford, 1996.

\item{\bf5.} A. Ivi\'c, `Large values of the error term in the
divisor problem', {\it Invent. Math.} {\bf71}(1983), 513-520.

\item{\bf6.} A. Ivi\'c, `The Riemann zeta-function'', John Wiley \&
Sons, New York, 1985.

\item{\bf7.} A. Ivi\'c, `On some problems involving the mean square
of $|\zt|$', {\it Bull.} CXVI {\it Acad. Serbe} 1998,
 {\it Classe des Sciences
math\'ematiques }{\bf23}, 71-76.

\item{\bf8.} Y.-K. Lau,  `Error terms in the summatory formulas
for certain number-theoretic functions', Doctoral dissertation,
University of Hong Kong, Hong Kong, 1999.

\item{\bf9.} T. Meurman, `On the mean square of the Riemann
zeta-function', {\it Quart. J. Math. Oxford} (2)38 (1987), 337-343.

\item{\bf10.}  W.G. Nowak, `On the divisor problem: $\D(x)$ over short
intervals', {\it Acta Arith.} {\bf109}(2003), 329-341.

\item{\bf11.} O. Robert and P. Sargos,  `Three-dimensional
exponential sums with monomials', {\it J. reine angew. Math.} (in print).

\item{\bf12.} E.C. Titchmarsh, `The theory of the Riemann zeta-function'
(2nd ed.),  University Press, Oxford, 1986.

\item{\bf13.} K.-M. Tsang, `Higher power moments of $\D(x)$, $E(t)$
and $P(x)$', {\it Proc. London Math. Soc. }(3){\bf 65}(1992), 65-84.

\item{\bf14.} W. Zhai, `On higher-power moments of $\D(x)$', and II.
{\it Acta Arith.} {\bf114}(2004), 35-54.

\endRefs

\bigskip
Aleksandar Ivi\'c

 Katedra Matematike RGF-a

Universitet u Beogradu

\DJ u\v sina 7, 11000 Beograd

Serbia and Montenegero, ivic\@rgf.bg.ac.yu

\bigskip

Patrick Sargos

 Institut Elie Cartan,

Universit\'e Henri Poincar\'e, BP. 239

54506 Nancy, France, sargos\@iecn.u-nancy.fr

\vfill


\bye